\documentclass[journal]{IEEEtran}

\usepackage{amssymb,amsmath}
\usepackage{color}
\usepackage{graphicx}
\usepackage{amscd}

\newtheorem{theorem}{Theorem}

\newcommand{\UC}{\mathcal{U}}%
\newcommand{\EC}{\mathcal{E}}%
\newcommand{\CC}{\mathcal{C}}%
\newcommand{\OC}{\mathcal{O}}%
\newcommand{\WC}{\mathcal{W}}%
\newcommand{\SC}{\mathcal{S}}%
\newcommand{\tm}{\times}%
\newcommand{\R}{\mathbb{R}}%
\newcommand{\rmd}{\mathrm{D}}%
\newcommand{\inv}{\mathrm{inv}}%
\newcommand{\cl}{\mathrm{cl}}%
\newcommand{\ep}{\varepsilon}%
\newcommand{\N}{\mathbb{N}}%
\newcommand{\rme}{\mathrm{e}}%
\newcommand{\inner}{\mathrm{int}}%
\renewcommand{\P}{\mathbb{P}}%
\newcommand{\F}{\mathbb{F}}%
\newcommand{\GL}{\mathrm{GL}}%
\newcommand{\SL}{\mathrm{SL}}%

\begin{document}

%
% paper title
% Titles are generally capitalized except for words such as a, an, and, as,
% at, but, by, for, in, nor, of, on, or, the, to and up, which are usually
% not capitalized unless they are the first or last word of the title.
% Linebreaks \\ can be used within to get better formatting as desired.
% Do not put math or special symbols in the title.
\title{Uniformly hyperbolic control theory}
%
%
% author names and IEEE memberships
% note positions of commas and nonbreaking spaces ( ~ ) LaTeX will not break
% a structure at a ~ so this keeps an author's name from being broken across
% two lines.
% use \thanks{} to gain access to the first footnote area
% a separate \thanks must be used for each paragraph as LaTeX2e's \thanks
% was not built to handle multiple paragraphs
%

\author{Christoph~Kawan
\thanks{The author thanks Professor Lai-Sang Young for the hospitality at the Courant Institute of Mathematical Sciences in 2014 and many inspiring mathematical conversations. Furthermore, he kindly acknowledges the support of DFG fellowship KA 3893/1-1 during this time. Moreover, the comments of the anonymous reviewers, which led to a substantial improvement of the paper, are appreciated very much. CK is with the Faculty of Computer Science and Mathematics, University of Passau, 94032 Passau, Germany; e-mail: christoph.kawan@uni-passau.de}}

% The paper headers
\markboth{Hyperbolic control theory}
{Shell \MakeLowercase{\textit{et al.}}: Bare Demo of IEEEtran.cls for IEEE Journals}
% The only time the second header will appear is for the odd numbered pages
% after the title page when using the twoside option.
% 
% *** Note that you probably will NOT want to include the author's ***
% *** name in the headers of peer review papers.                   ***
% You can use \ifCLASSOPTIONpeerreview for conditional compilation here if
% you desire.

% If you want to put a publisher's ID mark on the page you can do it like
% this:
%\IEEEpubid{0000--0000/00\$00.00~\copyright~2015 IEEE}
% Remember, if you use this you must call \IEEEpubidadjcol in the second
% column for its text to clear the IEEEpubid mark.

% use for special paper notices
%\IEEEspecialpapernotice{(Invited Paper)}

% make the title area
\maketitle

% As a general rule, do not put math, special symbols or citations
% in the abstract or keywords.
\begin{abstract}
This paper gives a summary of a body of work at the intersection of control theory and smooth nonlinear dynamics. The main idea is to transfer the concept of uniform hyperbolicity, central to the theory of smooth dynamical systems, to control-affine systems. Combining the strength of geometric control theory and the hyperbolic theory of dynamical systems, it is possible to deduce control-theoretic results of non-local nature that reveal remarkable analogies to the classical hyperbolic theory of dynamical systems. This includes results on controllability, robustness, and practical stabilizability in a networked control framework.%
\end{abstract}

% Note that keywords are not normally used for peerreview papers.

\iffalse
\fi

\section{Introduction}

The concept of uniform hyperbolicity, introduced in the 1960s by Stephen Smale, has become a cornerstone for the hyperbolic theory of dynamical systems, developed in the ensuing decades. This concept, which axiomatizes the geometric picture behind the horseshoe map and other complex systems, has been successfully generalized in various directions not long after its introduction to analyze a broad variety of systems (e.g., to non-uniform hyperbolicity, partial hyperbolicity and dominated splittings). A uniformly hyperbolic (discrete-time) system is essentially characterized by the fact that the linearization along any of its orbits behaves like a linear operator without eigenvalues on the unit circle, i.e., by a splitting of each tangent space into a direct sum of a stable and an unstable eigenspace. The uniformity is expressed by a uniform estimate on the contraction and expansion rates. We refer to \cite{Has} for a comprehensive survey of results related to hyperbolic dynamical systems.%

Uniform hyperbolicity and its generalizations also occur quite naturally in nonlinear control systems, which calls for a systematic transfer of the methods developed for the analysis of hyperbolic dynamical systems to control systems in order to gain new insights in control-theoretic problems. However, so far not much effort has been put into the development of a `hyperbolic control theory'. The aim of this paper is to provide a survey of the existing results, which show that a combination of techniques from geometric control theory and the uniformly hyperbolic theory of dynamical systems can lead to deep insights about global and semiglobal properties of control-affine systems with a compact and convex control range.%

These results are grounded on the topological theory of Colonius-Kliemann \cite{CKl} which provides an approach to understanding the global controllability structure of control systems. Two central notions of this theory are control and chain control sets. Control sets are the maximal regions of complete approximate controllability in the state space. The definition of chain control sets involves the concept of $\ep$-chains (also called $\ep$-pseudo-orbits) from the theory of dynamical systems. The main motivation for this concept comes from the facts that (i) chain control sets are outer approximations of control sets and (ii) chain control sets in general are easier to determine than control sets (both analytically and numerically).%

As examples show, chain control sets can support uniformly hyperbolic and, more generally, partially hyperbolic structures. For instance, every chain control set of a control-affine system on a flag manifold of a noncompact real semisimple Lie group, induced by a right-invariant system on the group, admits a partially hyperbolic structure, i.e., an invariant splitting of the tangent bundle into a stable, an unstable and a central subbundle. The paper \cite{DK2} provides a complete classification of those chain control sets on flags which are uniformly hyperbolic, using extensively the semigroup theory developed by San Martin and co-workers \cite{SM1,SM2,SM3,SM4}. Another way how a uniformly hyperbolic chain control set can arise is by adding sufficiently small control terms to an uncontrolled equation with a uniformly hyperbolic invariant set. In this case, under some control-theoretic regularity assumptions, the uniformly hyperbolic invariant set blows up to a uniformly hyperbolic chain control set.%

In the case of a uniformly hyperbolic chain control set, tools from the theory of smooth dynamical systems have been be applied to analyze controllability and robustness properties. In particular, it has been proved that a uniformly hyperbolic chain control set is the closure of a control set under the assumption of local accessibility, cf.~\cite{CDu}. As a consequence, complete controllability holds on the interior of the chain control set and the chain control set varies continuously in the Hausdorff metric in dependence on system parameters.%

Another control application of uniformly hyperbolic theory concerns the problem of practical stabilization under information constraints. Stabilization problems involving a communication channel of finite capacity which provides the controller with state information, have been studied by many authors (see, e.g., the survey \cite{Nea} and the monographs \cite{K1,MSa,YBa}). The main theoretical problem here is to determine the smallest capacity above which the stabilization objective can be achieved. For practical stabilization in the sense of rendering a compact subset $Q$ of the state space invariant, the notion of invariance entropy $h_{\inv}(Q)$ was introduced in \cite{CKa} as a measure for the associated critical channel capacity. This quantity measures the exponential complexity of the control task of keeping the system inside $Q$. In \cite{DK1} a formula for the invariance entropy $h_{\inv}(Q)$ of a uniformly hyperbolic chain control set $Q$ has been derived in terms of unstable volume growth rates along trajectories in $Q$. The proof of this formula in particular reveals the interesting fact that in order to make $Q$ invariant with a capacity arbitrarily close to $h_{\inv}(Q)$, control strategies that stabilize a periodic orbit in $Q$ are as good as any other strategy, thus this class of strategies is optimal.% 

The paper \cite{CLe} gives an application of this result to a problem related with a continuously stirred tank reactor. Moreover, the paper \cite{DK3} shows that the invariance entropy on uniformly hyperbolic chain control sets depends continuously on system parameters.%

In the following Sections \ref{sec_csets}--\ref{sec_tankreactor}, we explain these results in greater detail. In Section \ref{sec_coc}, we give a brief account of the related subjects known as `control of chaos' and `partial chaos', and in Section \ref{sec_futdir} we outline some problems and ideas for future research.%

{\bf Notation:} We write $\cl A$ and $\inner A$ for the closure and the interior of a set $A$, respectively. If $M$ is a smooth manifold, we write $T_xM$ for the tangent space to $M$ at $x$, and $TM$ for the tangent bundle of $M$. If $f:M\rightarrow N$ is a smooth map between manifolds, $\rmd f(x):T_xM \rightarrow T_{f(x)}N$ denotes its derivative at $x\in M$.%

\section{Control sets and chain control sets}\label{sec_csets}%

A control-affine system is governed by differential equations of the form%
\begin{equation}\label{eq_cas}
  \Sigma:\ \dot{x}(t) = f_0(x(t)) + \sum_{i=1}^m u_i(t)f_i(x(t)),\ u\in\UC,%
\end{equation}
where $x(t)$ lives on a Riemannian manifold $M$ (the state space) and $\UC$ is the set of admissible control functions, which we assume to be of the form $\UC=L^{\infty}(\R,U)$ with $U\subset\R^m$ being a compact and convex set satisfying $0\in\inner U$. Assuming that $f_0,f_1,\ldots,f_m$ are $\CC^1$-vector fields and that the unique solution $\varphi(t,x,u)$ for the initial value $x$ at time $t_0=0$ and the control $u$ exists for all $t\in\R$, regardless of $(u,x) \in \UC \tm M$, we obtain a skew-product flow (i.e., a flow of triangular structure)%
\begin{equation*}
  \Phi_t(u,x) = (\theta_tu,\varphi(t,x,u)),\ t\in\R,%
\end{equation*}
that acts on the extended state space $\UC \tm M$. Here%
\begin{equation*}
  \theta_tu = u(t+\cdot),\quad \theta_t:\UC \rightarrow \UC,\ t\in\R,%
\end{equation*}
denotes the shift flow on $\UC$. With the weak$^*$-topology of $L^{\infty}(\R,\R^m) = L^1(\R,\R^m)^*$, $\UC$ becomes a compact metrizable space and $\Phi$ a continuous flow, called the \emph{control flow} of $\Sigma$. We write $\varphi_{t,u} = \varphi(t,\cdot,u)$.%

A \emph{control set} of $\Sigma$ is a subset $D\subset M$ such that%
\begin{enumerate}
\item[(i)] for every $x\in D$ there is $u\in\UC$ with $\varphi(\R_+,x,u)\subset D$,%
\item[(ii)] for all $x,y\in D$ and every neighborhood $N$ of $y$ there are $u\in\UC$ and $T>0$ with $\varphi(T,x,u) \in N$ (i.e., approximate controllability holds on $D$), and%
\item[(iii)] $D$ is maximal with (i) and (ii) in the sense of set inclusion.%
\end{enumerate}
A \emph{chain control set} $E\subset M$ is a set such that%
\begin{enumerate}
\item[(i)] for every $x\in E$ there is $u\in\UC$ with $\varphi(\R,x,u)\subset E$,%
\item[(ii)] for all $x,y\in E$ and $\ep,T>0$ there are $n\in\N$, $u_0,\ldots,u_{n-1}\in\UC$, $x = x_0,x_1,\ldots,x_n = y$ and $t_0,t_1,\ldots,t_{n-1}\geq T$ such that $d(\varphi(t_i,x_i,u_i),x_{i+1}) < \ep$ for $i=0,1,\ldots,n-1$, and%
\item[(iii)] $E$ is maximal with (i) and (ii) in the sense of set inclusion.%
\end{enumerate}

Before we proceed, for the convenience of the reader, we explain the concept of chain transitivity used in the topological theory of dynamical systems to analyze recurrence properties (see also \cite[App.~B]{CKl}). If $\phi:\R \tm X \rightarrow X$ is a continuous flow on a metric space $(X,d)$, a set $A\subset X$ is called \emph{chain transitive} if for all $x,y\in A$ and $\ep,T>0$ there exists an \emph{$(\ep,T)$-chain} from $x$ to $y$, i.e., there are $n\in\N$, points $x = x_0,x_1,\ldots,x_n=y$, and times $t_0,t_1,\ldots,t_{n-1}\geq T$ so that $d(\phi(t_i,x_i),x_{i+1}) < \ep$ for $i=0,\ldots,n-1$. A point $x\in X$ is called \emph{chain recurrent} if for all $\ep,T>0$ there is an $(\ep,T)$-chain from $x$ to $x$. If $X$ is compact, then the set $R(\phi)$ of all chain recurrent points is closed and invariant. Moreover, the connected components of $R(\phi)$ are precisely the maximal invariant chain transitive sets and are called \emph{chain recurrent components}. The chain recurrent set essentially contains all relevant dynamical information of the flow. For instance, all $\alpha$- and $\omega$-limit set are contained in $R(\phi)$.%

Now we consider again the control-affine system \eqref{eq_cas}. The \emph{lift} of a chain control set $E$ is defined by%
\begin{equation*}
  \EC := \left\{(u,x)\in\UC\tm M\ :\ \varphi(\R,x,u) \subset E\right\}.%
\end{equation*}
It is a maximal invariant chain transitive set of the control flow, hence a chain recurrent component if $M$ is compact. If $\Sigma$ is locally accessible and $D$ is a control set with nonempty interior, then $D$ is contained in a chain control set (which is unique, since different chain control sets are disjoint). In general, chain control sets are closed, while control sets are neither open nor closed except when they are invariant in backward or forward time, respectively.%

A chain control set $E$ is \emph{uniformly hyperbolic without center bundle} if it is compact and for every $(u,x)\in\EC$ there exists a splitting%
\begin{equation*}
  T_xM = E^-_{u,x} \oplus E^+_{u,x}%
\end{equation*}
into linear subspaces such that%
\begin{enumerate}
\item[(i)] $\rmd\varphi_{t,u}(x)E^{\pm}_{u,x} = E^{\pm}_{\Phi_t(u,x)}$ for all $t\in\R$ and $(u,x)\in\EC$, and%
\item[(ii)] there are constants $c,\lambda>0$ such that%
\begin{equation*}
  |\rmd\varphi_{t,u}(x)v| \leq c^{-1}\rme^{-\lambda t}|v| \mbox{\quad for all\ } t\geq 0,\ v\in E^-_{u,x}%
\end{equation*}
and%
\begin{equation*}
  |\rmd\varphi_{t,u}(x)v| \geq c\rme^{\lambda t}|v| \mbox{\quad for all\ } t\geq 0,\ v\in E^+_{u,x}.%
\end{equation*}
\end{enumerate}
This definition is independent of the Riemannian metric, however, the constant $c$ depends on the choice of the metric. From the two conditions it automatically follows that the subspaces $E^{\pm}_{u,x}$ change continuously with $(u,x)$, cf.~\cite[Ch.~6]{K1}.%

If the tangent spaces $T_xM$, $(u,x)\in\EC$, admit continuous invariant splittings%
\begin{equation*}
  T_xM = E^-_{u,x} \oplus E^0_{u,x} \oplus E^+_{u,x},%
\end{equation*}
with uniform exponential contraction on $E^-_{u,x}$ and expansion on $E^+_{u,x}$ (as above), and $E^0_{u,x}$ is one-dimensional and corresponds to the flow direction for constant controls, we say that $E$ is \emph{uniformly hyperbolic with center bundle}.%

These are natural extensions of the well-known concept of a uniformly hyperbolic invariant set in dynamical systems, cf.~\cite{Has}. Similar extensions for skew-product systems have been studied before in the context of random dynamical systems \cite{GKi,Liu} and almost periodic differential equations \cite{MZh}.%

In the rest of the paper, we will often use the abbreviation \emph{u.h.}\ for \emph{uniformly hyperbolic}.%

\section{A structural result}\label{sec_sr}%

For uncontrolled time-invariant systems in continuous time, the notion of uniform hyperbolicity without center bundle has not much meaning, because any trajectory that is bounded and bounded away from equilibria allows for neither exponential contraction nor expansion in the flow direction. As a consequence, every uniformly hyperbolic invariant set without center bundle is a discrete set of equilibrium points (a precise proof for this well-known fact can be found in \cite[Prop.~4]{K2}). However, for control-affine systems the situation is different, since uniformly hyperbolic chain control sets with nonempty interior exist. However, their structure can be shown to be relatively simple; under mild assumptions, their lifts are graphs over $\UC$. More precisely, the following theorem holds, cf.~\cite[Thm.~5]{K2}.%

\begin{theorem}\label{thm_str}
Let $E$ be a u.h.~chain control set without center bundle. Assume that $\EC$ is an isolated invariant set of $\Phi$ and let $u_0$ be a constant control function with value in $\inner U$. Additionally suppose that the following hypotheses are satisfied:%
\begin{enumerate}
\item[(i)] The vector fields $f_0,f_1,\ldots,f_m$ are of class $\CC^{\infty}$ and the Lie algebra generated by them has full rank at each point of $E$.%
\item[(ii)] For each $x$ with $(u_0,x)\in\EC$ and each $\rho\in(0,1]$ it holds that $x\in\inner\OC^+_{\rho}(x)$, where $\OC^+_{\rho}(x)=\{\varphi(t,x,u) : t\geq 0,\ u\in\UC^{\rho}\}$ with
\begin{equation*}
  \UC^{\rho} = \{u \in \UC\ :\ u(t) \in u_0 + \rho(U - u_0) \mbox{ a.e.}\}.%
\end{equation*}
\end{enumerate}
Then $\EC$ is the graph of a continuous function $\UC \rightarrow E$.%
\end{theorem}

The condition that $\EC$ is isolated invariant means that $\EC$ is the largest compact invariant set in a neighborhood of $\EC$.%

For condition (ii) note that each $x\in E$ kept in $E$ by a constant control function is necessarily an equilibrium point. Since hyperbolic equilibria are isolated, condition (ii) has to be checked only for finitely many points $x$. A sufficient condition, independent of $\rho$, for (ii) to hold is the controllability of the linearization at $(u_0,x)$.%

The above theorem can be seen as a technical lemma, which is very useful for proving more advanced results such as Theorem \ref{thm_ie} in Section \ref{sec_ie}, characterizing the smallest bit rate in a digital channel above which $E$ can be rendered invariant. Theorem \ref{thm_str} implies that the restriction of the control flow $\Phi$ to $\EC$ is topologically conjugate to the shift flow on $\UC$ via the continuous projection map $\pi_1:(u,x) \mapsto u$, i.e., $\pi_1:\EC\rightarrow\UC$ is a homeomorphism and the following diagram commutes:%
\begin{equation*}
  \begin{CD}
    \EC @>\Phi_t>>\EC\\
    @V \pi_1 VV @VV \pi_1 V\\
     \UC @>>\theta_t> \UC
  \end{CD}
\end{equation*}
Indeed, the characterization of $\EC$ as a graph $\{(u,x(u)) : u \in \UC\}$ implies that $\pi_1:\EC\rightarrow\UC$ is invertible with $\pi_1^{-1}(u) = (u,x(u))$. Since both $\EC$ and $\UC$ are compact metric spaces, $\pi_1$ is a homeomorphism. Moreover, $\pi_1(\Phi_t(u,x(u))) = \theta_tu = \theta_t \pi_1(u,x(u))$, showing that $\pi_1$ is a conjugacy. Hence, all topological properties of the shift flow carry over to $\Phi|_{\EC}$.%
 
\section{Algebraic examples}\label{sec_ae}%

Examples of u.h.~chain control sets without center bundle can be constructed as follows. We start with a bilinear control system on $\R^{n+1}$:%
\begin{equation}\label{eq_bis}
  \dot{x}(t) = \Bigl(A_0 + \sum_{i=1}^m u_i(t)A_i\Bigr)x(t),\ u\in\UC.%
\end{equation}
Since $\varphi_{t,u}$ is a linear isomorphism of $\R^{n+1}$ for all $t$ and $u$, the system induces another control-affine system on the $n$-dimensional projective space $\P^n = \P(\R^{n+1})$, the space of all lines through $0$ in $\R^{n+1}$. A description of the chain control sets of this system comes out of Selgrade's theorem about linear flows on vector bundles with chain transitive base (cf.~\cite[Ch.~5]{CKl}). In our case, the vector bundle is $\UC \tm \R^{n+1}$ and the linear flow is the control flow of \eqref{eq_bis}. By chain transitivity of $\theta$, Selgrade's result can be applied and it yields a description of the chain recurrent components of the control flow on $\UC \tm \P^n$, whose projections to $M$ are the chain control sets. More precisely, there exists a Whitney sum decomposition%
\begin{equation*}
  \UC \tm \R^{n+1} = \WC^1 \oplus \cdots \oplus \WC^r,%
\end{equation*}
where $r \leq n+1$, into $\Phi$-invariant subbundles $\WC^i$ and the chain recurrent components $\EC_i$ correspond to these $\WC^i$ in the sense that%
\begin{equation*}
  \EC_i = \left\{ (u,\P x) \in \UC \tm \P^n : (u,x) \in \WC^i \right\},\ 1 \leq i \leq r.%
\end{equation*}
We fix $i\in\{1,\ldots,r\}$ and order $\WC^1,\ldots,\WC^r$ by increasing growth rates. Then we define%
\begin{align*}
  \WC^- &:= \WC^1 \oplus \cdots \oplus \WC^{i-1},\ \WC^+ := \WC^{i+1} \oplus \cdots \oplus \WC^r,\\
	\WC^0 &:= \WC^i.%
\end{align*}
Projecting to $\P^n$ by $E^{\pm}_{u,\P x} := \rmd\P(x)\WC^{\pm}(u)$ and $E^0_{u,\P x} := \rmd\P(x)\WC^0(u)$ (here $\P:\R^{n+1}\backslash\{0\}\rightarrow\P^n$ denotes the projection map which sends $x$ to its equivalence class, i.e., the line through $x$), we obtain a splitting%
\begin{equation}\label{eq_phs}
  T_{\P x}\P^n = E^-_{u,\P x} \oplus E^0_{u,\P x} \oplus E^+_{u,\P x}.%
\end{equation}
From the fact that $\WC^+$, $\WC^0$ and $\WC^-$ are exponentially separated one can deduce that the splitting $E^- \oplus E^+$ is uniformly hyperbolic. Hence, if $E^0$ is trivial, i.e., $0$-dimensional, the chain control set $E_i = \{\P x\in\P^n: (u,x)\in \EC_i\}$ turns out to be uniformly hyperbolic without center bundle.%

Chain control sets on $\P^n$ of this type can also be seen as control-dependent `eigenvectors'. Indeed, if $A_0$ has a real eigenvalue of multiplicity one and $x_0 \in \R^{n+1}$ is an associated eigenvector, then $\P x_0$ is an isolated hyperbolic equilibrium of the induced flow on $\P^n$. By adding small control terms $u_i A_i$, this equilibrium can blow up to a u.h.\ chain control set $E$ without center bundle with nonempty interior, and each point in $E$ is of the form $z(u)$ for an equivariant continuous map $z:\UC \rightarrow E$, where $z(0) = \P x_0$.%

A bilinear system not only induces a system on projective space, but also on Grassmannians and flag manifolds. Another way to describe these systems is by looking at a right-invariant system on the Lie group $G = \GL(n+1,\R)$ or $G = \SL(n+1,\R)$ (if $\mathrm{tr} A_i = 0$ for $i = 0,1,\ldots,m$) given by%
\begin{equation*}
  \dot{g}(t) = \Bigl(A_0 + \sum_{i=1}^m u_i(t)A_i\Bigr)g(t),\quad u\in\UC,%
\end{equation*}
where $g(t) \in G$, and viewing the flag manifolds as homogeneous spaces of this group. This can be generalized by replacing the group $\SL(n+1,\R)$ with an arbitrary non-compact semisimple Lie group $G$ and the bilinear system with a right-invariant system on $G$. Then one can study the induced systems on the generalized flag manifolds $\F_{\Theta} = G/P_{\Theta}$, where $P_{\Theta}$ denotes a parabolic subgroup of $G$ characterized by a set $\Theta$ of simple roots. With tools from semisimple Lie theory and semigroup theory one can show that the chain control sets of such systems also have a partially hyperbolic structure as in \eqref{eq_phs}. In \cite{DK2} the u.h.\ chain control sets without center bundle are characterized via the so-called flag type of the control flow. We do not go into further details here, because this would necessitate to introduce plenty of Lie-theoretic notions.%

We note that also in the classical uniformly hyperbolic theory the simplest examples are given by algebraic systems, namely by linear automorphisms of the $n$-dimensional torus, which are Anosov diffeomorphisms, i.e., diffeomorphisms that admit a hyperbolic structure on the whole state space, cf.~\cite{Has}.%

\section{Controllability and robustness}\label{sec_contr_rob}%

A very useful feature of uniformly hyperbolic systems is the shadowing property. Roughly speaking, shadowing means that $\ep$-close to any $\delta$-approximate orbit (as used in the definition of chain control sets) there exists a unique real orbit, where $\delta = \delta(\ep)$.%

Using the shadowing property, it is possible to show that a uniformly hyperbolic chain control set $E$ (with or without center bundle) is the closure of a control set, provided that it a has nonempty interior and local accessibility holds. In general, i.e., without assuming uniform hyperbolicity, this does not hold. In fact, a chain control set may contain several control sets that have positive distance to each other.%

This result has the remarkable consequence that compact uniformly hyperbolic chain control sets change continuously in the Hausdorff metric with respect to parameters of the systems. Such a parameter, e.g., could be the size of the control range or any vector parameter which smoothly influences the vector fields $f_0,f_1,\ldots,f_m$. More precisely, the following result holds, cf.~\cite[Thm.~2 and Thm.~3]{CDu}.%

\begin{theorem}\label{thm_robustness1}
Consider a parametrized control-affine system of the form%
\begin{equation*}
  \dot{x}(t) = f_0(\alpha,x(t)) + \sum_{i=1}^m u_i(t)f_i(\alpha,x(t)),\quad u\in\UC,%
\end{equation*}
where the parametrized vector fields $f_i:A\tm M \rightarrow TM$, with $A \subset \R^k$ and $M$ a smooth manifold, are of class $\CC^{\infty}$. Assume that for a fixed parameter $\alpha^0 \in \inner A$ the Lie algebra rank condition holds on a uniformly hyperbolic chain control set $E^{\alpha^0}$ (with or without center bundle). Then $E^{\alpha^0}$ is the closure of a control set and for each $\alpha$ in a neighborhood of $\alpha^0$ there is a unique control set $D^{\alpha}$ such that $\alpha\mapsto \cl D^{\alpha}$ is continuous in the Hausdorff metric at $\alpha^0$ with $\inner D^{\alpha}\cap \inner D^{\alpha^0} \neq \emptyset$.%
\end{theorem}

Note that this theorem in particular implies that a control-affine system is completely controllable on the interior of a uniformly hyperbolic chain control set $E$ provided that the Lie algebra rank condition is satisfied on $E$. This follows by combining the approximate controllability on the control set and the local accessibility guaranteed by the Lie algebra rank condition. Since chain control sets are easier to determine than control sets (numerically and analytically), this result is important for the understanding of the controllability properties of a system.%

Theorem \ref{thm_robustness1} corresponds to well-known robustness results in the uniformly hyperbolic theory of smooth dynamical systems (e.g., structural stability of Axiom A systems). The controllability part of the result is analogous to topological transitivity on Axiom A basic sets, cf.~\cite[Sec.~3.3(a)]{Has}.%

\section{Networked control and invariance entropy}\label{sec_ie}%

Networked control systems are spatially distributed systems, in which the communication between sensors, controllers and actuators is accomplished through a shared digital communication network. Examples can be found in vehicle tracking, underwater communications for remotely controlled surveillance and rescue submarines, remote surgery, space exploration and aircraft design. Another large field of applications can be found in modern industrial systems, where industrial production is combined with information and communication technology (`Industry 4.0').%

In networked control systems, the analog system outputs must be encoded in finite bit strings to be transmitted over the communication network. Realistic models of such systems therefore challenge the standard assumption of control theory that controllers and actuators have access to continuous-valued state information, i.e., information of infinite precision. As a consequence, the characterization of a property such as stabilizability for a networked control system involves not only characteristics of the dynamical system, but also of the communication network. In the simplest setup -- one dynamical system connected via a noiseless digital channel to a controller -- this reduces to the computation of the smallest channel capacity, above which the system can be stabilized.%

For practical stabilization (in the sense of set-invariance), the invariance entropy provides a measure for this critical channel capacity. Formally, the invariance entropy is a non-negative quantity that can be assigned to any compact controlled invariant set $Q$ and measures the exponential complexity of the control task of keeping the system inside $Q$, cf.~\cite{CKa,K1}. For continuous-time systems, it is defined as follows. Let $K$ be a compact subset of $Q$. For $\tau>0$, a set $\SC\subset\UC$ of controls is called $(\tau,K,Q)$-spanning if for each $x\in K$ there is $u\in\SC$ with $\varphi([0,\tau],x,u) \subset Q$. The minimal cardinality of such a set is denoted by $r_{\inv}(\tau,K,Q)$ and the exponential growth rate%
\begin{equation*}
  h_{\inv}(K,Q) := \limsup_{\tau\rightarrow\infty}\frac{1}{\tau}\log r_{\inv}(\tau,K,Q)%
\end{equation*}
is called the \emph{invariance entropy} of the pair $(K,Q)$. Invariance entropy is essentially equivalent to the notion of \emph{topological feedback entropy}, introduced in \cite{Nai}.%

\begin{figure}[hp]
	\begin{center}
		\includegraphics[width=8.10cm,height=2.40cm,angle=0]{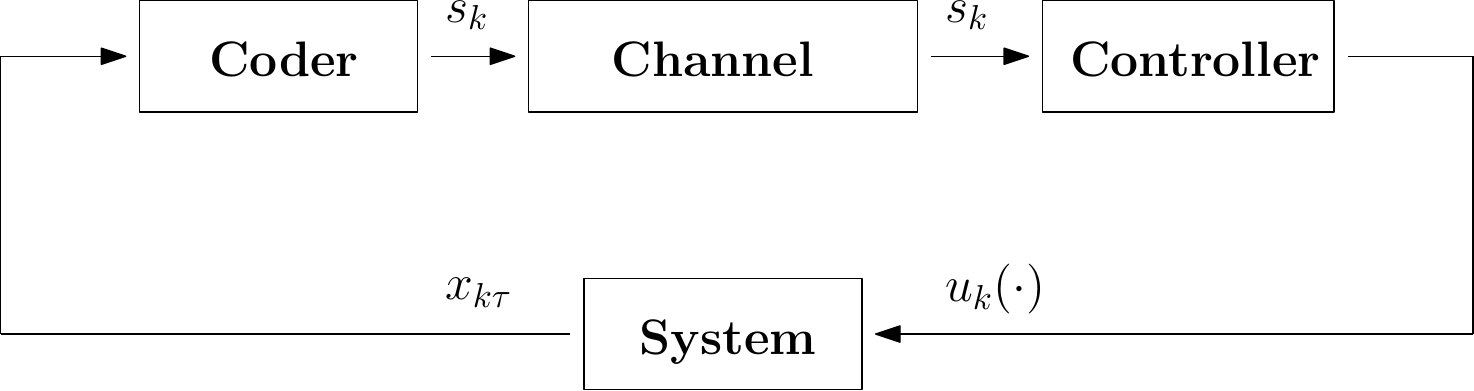}
	\end{center}
\caption{Control over a digital channel}\label{fig1}
\end{figure}

The information-theoretic interpretation of $h_{\inv}(K,Q)$, mentioned above, can be explained as follows. Suppose that a sensor measures the states of the system at discrete sampling times $\tau_k = k \tau$ for some $\tau>0$. A coder receiving these measurements generates at each sampling time $\tau_k$ a symbol $s_k$ from a finite coding alphabet $S_k$ of time-varying size. This symbol is transmitted through a digital noiseless channel to a controller, see Fig.~\ref{fig1}. The controller, upon receiving $s_k$, generates an open-loop control $u_k$ on $[0,\tau]$ used as the control input in the time interval $[\tau_k,\tau_{k+1}]$. The aim of this coding and control device is to keep trajectories in $Q$, when $x_0 \in K$. The invariance entropy $h_{\inv}(K,Q)$ is the smallest channel bit rate above which a coder and a controller can be designed so that this control objective is achieved.%

For u.h.~chain control sets without center bundle we have the following result, cf.~\cite[Thm.~5.4]{DK1}.%

\begin{theorem}\label{thm_ie}
Assume that $E$ is a u.h.\ chain control set without center bundle of $\Sigma$ with nonempty interior satisfying the conditions of Theorem \ref{thm_str}. Then $E$ is the closure of a control set $D$ and for every compact $K\subset D$ with positive volume it holds that%
\begin{equation}\label{eq_ief}
  h_{\inv}(K,E) = \inf_{(u,x)\in\EC}\limsup_{\tau\rightarrow\infty}\frac{1}{\tau}\log J^+_x\rmd\varphi_{\tau,u},%
\end{equation}
where%
\begin{equation*}
  J^+_x\rmd\varphi_{\tau,u} = \left|\det\left(\rmd\varphi_{\tau,u}(x)|_{E^+_{u,x}}:E^+_{u,x}\rightarrow E^+_{\Phi_{\tau}(u,x)}\right)\right|.%
\end{equation*}
\end{theorem}

Moreover, the infimum in \eqref{eq_ief} can be taken only over the $\Phi$-periodic points $(u,x)\in\EC$. For a periodic point $(u,x)$ the $\limsup$ in \eqref{eq_ief} is equal to the sum of the positive Lyapunov exponents along the corresponding periodic trajectory. This sum in turn is the exponential complexity of the control task of stabilizing the system (locally and exponentially) around the periodic trajectory.%

Again, we have a remarkable analogy to the theory of smooth dynamical systems, where the notions of measure-theoretic and topological entropy play a crucial role, and can be shown to depend only on the periodic orbits in the case of uniform hyperbolicity, cf.~\cite{Bo1,Bo2}.%

Additionally to the robustness result of Theorem \ref{thm_robustness1}, concerning the dependence of u.h.\ chain control sets on parameters, a robustness result for the associated invariance entropy can be proved, which reads as follows, cf.~\cite{DK3}.%

\begin{theorem}
Consider a parametrized control-affine system of the form%
\begin{equation*}
  \Sigma^{\alpha}:\ \dot{x}(t) = f_0(\alpha,x(t)) + \sum_{i=1}^m u_i(t)f_i(\alpha,x(t)),\quad u\in\UC,%
\end{equation*}
where the parametrized vector fields $f_i:A\tm M \rightarrow TM$, with $A \subset \R^k$ and $M$ a smooth manifold, are of class $\CC^{\infty}$. Assume that for a fixed parameter $\alpha^0 \in \inner A$ a u.h.\ chain control set $E^{\alpha^0}$ without center bundle exists, which satisfies the assumptions of Theorem \ref{thm_str}. Furthermore, assume that $\Sigma^{\alpha^0}$ satisfies the Lie algebra rank condition on $E^{\alpha^0}$. Then for each $\alpha$ in a neighborhood of $\alpha^0$, the system $\Sigma^{\alpha}$ has a u.h.\ chain control set $E^{\alpha}$ without center bundle and the map $\alpha \mapsto h_{\inv}(K,E^{\alpha})$ is continuous on this neighborhood.%
\end{theorem}

\section{Application to a stirred tank reactor}\label{sec_tankreactor}

In this section, we provide an example which comes from a concrete application. The equation%
\begin{align*}
  \left(\begin{array}{c} \dot{x}_1 \\ \dot{x}_2 \end{array}\right) &= \left(\begin{array}{c} -x_1 - a(x_1 - x_c) + B\alpha(1 - x_2)\rme^{x_1}\\ -x_2 + \alpha(1-x_2)\rme^{x_1} \end{array}\right)\\
	& \qquad\qquad\qquad + u(t)\left(\begin{array}{c} x_c - x_1 \\ 0 \end{array}\right)%
\end{align*}
models a continuous stirred tank reactor with Arrhenius' dynamics, cf., e.g.~\cite{Poo}.%

\begin{figure}[hp]
	\begin{center}
		\includegraphics[width=6.50cm,height=5.50cm,angle=0]{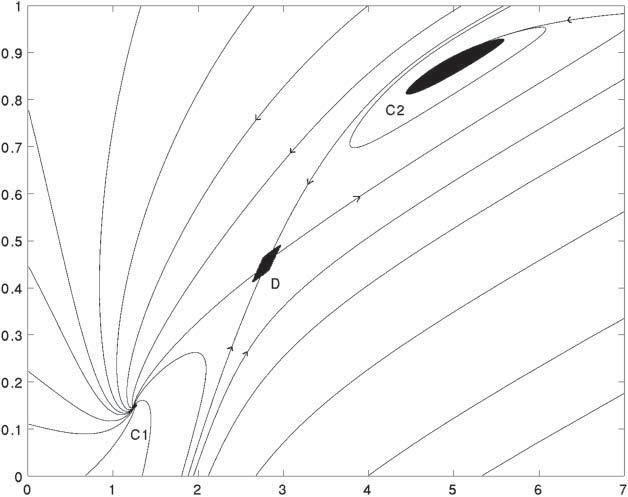}
	\end{center}
\caption{Phase portrait of the continuous flow stirred tank reactor and control sets}\label{fig2}
\end{figure}

Here $x_1$ is the (dimensionless) temperature; $x_2$ is the product concentration; and $a,\alpha,B,x_c$ are positive constants. The parameter $x_c$ is the coolant temperature, and hence the control affects the heat transfer coefficient. In the following, we look at the system for the parameters%
\begin{align*}
  a &= 0.15,\ \alpha = 0.05,\ B = 7.0,\ x_c = 1.0,\\
	&\qquad U_{\rho} = [−\rho,\rho] \mbox{ with } 0 < \rho \leq 0.15.%
\end{align*}
Because of the physical constraints, we consider the system in the set $[0,\infty) \tm [0,1] \subset \R^2$. For each fixed $u\in U_{\rho}$ (i.e., for each constant control function), we have three equilibria, given by%
\begin{equation*}
  p_i = (z_i,y_i),\quad y_i = \frac{\alpha \rme^{z_i}}{1 + \alpha z_i}%
\end{equation*}
for $i=0,1,2$, where $z_1 < z_0 < z_2$ are the solutions of the transcendental equation%
\begin{equation*}
  -z - (a+u)(z-x_c) + B\alpha\frac{\rme^z}{1 + \alpha\rme^z} = 0.%
\end{equation*}
The equilibria $p_1$ and $p_2$ are stable, while $p_0$ is hyperbolic, i.e., the linearization at $p_0$ has one positive and one negative eigenvalue. The phase portrait of the uncontrolled equation is depicted in Fig.~\ref{fig2}.

The system satisfies the Lie algebra rank condition at every point of the forward-invariant set $(0,\infty) \tm (0,1)$ (see \cite[Sec.~9.1]{CKl} for a verification of this fact). Numerical computations suggest that the rectangle $[0,7] \tm [0,1]$ contains exactly three control sets $C_1^{\rho}$, $C_2^{\rho}$ and $D^{\rho}$, containing the equilibria $p_1(u),p_2(u)$ and $p_0(u)$ for $u \in \inner U^{\rho} = (-\rho,\rho)$ in their interiors (this follows from an application of \cite[Cor.~4.1.12]{CKl}). The control sets $C_1^{\rho}$ and $C_2^{\rho}$ are invariant, while $D^{\rho}$ is variant (i.e., escape from $D^{\rho}$ is possible). Figure \ref{fig2} shows the situation for $\rho = 0.15$.%

An interesting property of this system is that the stable equilibrium $p_2$ with the highest product concentration cannot be realized for technical reasons (see \cite{BBM}). Hence, it is of interest to keep the system near the hyperbolic equilibrium point $p_0$. For $\rho$ small enough, the control set $D^{\rho}$ is uniformly hyperbolic, which follows from standard results on the persistence of hyperbolicity under small perturbations (see, e.g., \cite{Liu}). From \cite[Cor.~3.4.10]{CKl} it follows that the chain control sets of the system shrink to the equilibrium points $p_i$ as the control range shrinks to $\{0\}$. Hence, for $\rho$ small enough, the chain control set $E_0^{\rho}$ (containing $p_0$) is uniformly hyperbolic without center bundle and thus satisfies $E_0^{\rho} = \cl D^{\rho}$. We can therefore apply all the results of the preceding sections, which in particular tell us the following about the system under consideration:%
\begin{enumerate}
\item[(1)] The control set $D^{\rho}$ varies continuously in the Hausdorff metric when the parameters $(\rho,a,\alpha,B,x_c)$ are varied.% 
\item[(2)] The critical bit rate necessary for rendering $D^{\rho}$ invariant varies continuously with the parameters.%  
\item[(3)] Rendering $D^{\rho}$ invariant with a bit rate arbitrarily close to the theoretical infimum is possible by stabilization of a periodic orbit in $\inner D^{\rho}$. The corresponding bit rate for this control strategy is given by the sum of the unstable eigenvalues of the monodromy operator associated with the periodic orbit (divided by the period).%
\end{enumerate}

\section{Control and partial control of chaos}\label{sec_coc}%

This survey on uniformly hyperbolic control theory would not be complete without some remarks about the subjects known by the names `control of chaos' and `partial control (of chaos)', which also use ideas from the hyperbolic theory of dynamical systems. The idea of the first is to produce desired controlled trajectories with low energy use by fixing an unstable periodic orbit inside a chaotic (possibly hyperbolic) attractor of an uncontrolled dynamical system and stabilizing this orbit via very small time-dependent perturbations of a system parameter, applied once in a while. This method (which essentially exists in two variants, called the \emph{OGY method} \cite{OGY} and the \emph{Pyragas method} \cite{Pyr}) has proven to be extremely useful and effective. Its main advantage is that it does not require a detailed model of the chaotic attractor, but only some information about a Poincar\'e section, which is used to determine the periodic orbit. Since a typical trajectory in a neighborhood of the chaotic attractor will come close to any periodic orbit, one can just wait until the system runs into a small enough neighborhood of the fixed trajectory and then apply the control algorithm. Experimental applications of control of chaos include turbulent fluids, oscillating chemical reactions, magneto-mechanical oscillators, and cardiac tissues.%

While the idea of controlling chaos was first introduced in the 1990 paper \cite{OGY}, the idea of `partial control' is relatively new, cf.~\cite{ZSY}. Here one fixes a non-attracting hyperbolic invariant set of an uncontrolled equation and adds both noise and control terms to the equation. Using horseshoe structures, it is possible to show the existence of a non-trivial `safe region' around the hyperbolic set which is controlled invariant in the sense that whenever the initial condition lies in this region, no matter what the disturbance is, one can always choose a control that keeps the system in the safe region for the next time step. The safe region, in general, has a complicated geometric structure, although it is not a fractal set. The main novelty in this approach is that the absolute value of the control can be chosen strictly smaller than the maximal possible disturbance. The name `partial control' comes from the fact that the method does not allow to follow any particular trajectory, but only guarantees that the system stays close to the chaotic set.%

While the theory described in Sections \ref{sec_csets}--\ref{sec_ie} is rather concerned with the general analysis of control systems, the ideas described in this section aim at concrete control algorithms to be practically implemented, even without knowing the exact equations governing the system.%

\section{Conclusion and future directions}\label{sec_futdir}%

In this paper, we reviewed a number of results about non-local controllability and robustness properties of control-affine systems. All of these results rely on a combination of methods from geometric control and techniques from the uniformly hyperbolic theory of dynamical systems.%

The described results mainly deal with the simpler type of uniform hyperbolicity without center bundle. It is desirable to extend the results to the discrete-time case and uniformly hyperbolic chain control sets with center bundle. Particular difficulties will arise in extending the result about invariance entropy, because its proof is heavily based on the structural result of Theorem \ref{thm_str}, which does no longer hold in the case with center bundle. Furthermore, the proof of Theorem \ref{thm_ie} requires an elaborate result of Coron \cite{Cor} on the genericity of universally regular control functions in smooth systems, an analogue of which does not seem to be available in the discrete-time case.%

Another possible future direction concerns the extension of the controllability and robustness result described in Section \ref{sec_contr_rob} to the partially hyperbolic case, when a center bundle of dimension $\geq 2$ is present. Examples for this behavior can be derived from the construction described in Section \ref{sec_ae}, or the more general Lie-theoretic construction in \cite{DK2}. In the uncontrolled case, the theory of partial hyperbolicity (see \cite{HPe}) is quite well-developed so that one can hope to transfer methods of this theory to the control-affine case and obtain new results about the controllability structure inside a partially hyperbolic chain control set.%

\begin{IEEEbiography}{Christoph Kawan}
Christoph Kawan received the diploma and doctoral degree in mathematics (both under supervision of Prof. Fritz Colonius) 
from the University of Augsburg, Germany, in 2006 and 2009, respectively. He did research at the State University of 
Campinas, Brazil, in 2011. In 2014, he spent nine months as a research scholar at the Courant Institute of Mathematical 
Sciences at New York University. He is the author of the book `Invariance Entropy for Deterministic Control Systems - An 
Introduction' (Lecture Notes in Mathematics 2089. Berlin: Springer, 2013). Since January 2015, he works in the research 
group of Prof. Fabian Wirth at the University of Passau, Germany. His research interests include networked control, the 
qualitative theory of nonautonomous dynamical systems and the topological classification of dynamical systems.
\end{IEEEbiography}

\end{document}